\documentclass[11pt]{amsart}   
\usepackage{amssymb, amsfonts, amsmath, xcolor}
\usepackage[margin=1 in]{geometry}
\newcommand{\CC}{\mathbb{C}}

\newcommand{\RR}{\mathbb{R}}

\newcommand{\stern}[2]{\mathchoice
{\left\langle \begin{matrix} #1 \\ #2 \end{matrix} \right\rangle}
{\left\langle \begin{smallmatrix} #1 \\ #2 \end{smallmatrix} \right\rangle}
{\left\langle \begin{smallmatrix} #1 \\ #2 \end{smallmatrix} \right\rangle}
{\left\langle \begin{smallmatrix} #1 \\ #2 \end{smallmatrix} \right\rangle}
}

\newenvironment{sbm} 
    {\left[ \begin{smallmatrix}
    }
    { 
     \end{smallmatrix} \right]
    }

\newcommand{\sym}{\mathrm{sym}}

\newtheorem*{eg}{Example}
\newtheorem*{MT}{Main Theorem}
\newtheorem*{MC}{Consequences of the Main Theorem}

\begin{document}

\title{Proof of a conjecture of Stanley about Stern's array}
\author{David E Speyer}
\maketitle

\begin{abstract}
Stanley, building on work of Stern, defined an array of numbers by the recurrence $\stern{n}{2k} = \stern{n-1}{k}$, $\stern{n}{2k+1} = \stern{n-1}{k} + \stern{n-1}{k+1}$. 
Stanley showed that, for each positive integer $r$, the sequence $s_n^r:= \sum_k \stern{n}{k}\!^r$ obeys a homogeneous linear recurrence in $n$ of length $r/2+O(1)$. 
Numerical evidence, however, suggested that $s_n^r$ obeys shorter recurrences, of length  $r/3+O(1)$. 
We prove Stanley's conjecture.
\end{abstract}

Richard Stanley~\cite{Stan}, following a construction of Stern~\cite{Stern}, introduced an array of numbers defined as follows:
We start with the sequence $\cdots 0001000 \cdots$. 
We copy down this sequence and, in between each two consecutive elements of this sequence, we insert the sum of those two elements to obtain
$\cdots 0 \textcolor{red}{0} 0  \textcolor{red}{0} 0  \textcolor{red}{1} 1 \textcolor{red}{1}0 \textcolor{red}{0} 0  \textcolor{red}{0} 0 \cdots$; the new elements are shown in red. We repeat this process to obtain the array shown below (with $0$'s omitted).
\[ \begin{array}{ccccccccccccccc}
&&&&&&&1&&&&&&& \\
&&&1&&&&1&&&&1&&& \\
&1&&1&&2&&1&&2&&1&&1& \\
1&1&2&1&3&2&3&1&3&2&3&1&2&1&1 \\
\end{array} \]
We label the nonzero elements of the $n$-th row as $\stern{n}{k}$ for $1 \leq k \leq 2^n-1$. So we have the recurrence:
\[ \stern{n}{2k} = \stern{n-1}{k}  \ \mbox{and}\   \stern{n}{2k+1} = \stern{n-1}{k} + \stern{n-1}{k+1} \]
with initial conditions
\[ \stern{1}{1} = 1 \ \mbox{and}\   \stern{1}{k}=0 \ \mbox{for}\ k \neq 1 . \]

For binomial coefficients, we have the easy identities $\sum_k \binom{n}{k} = 2^n$ and $\sum_k \binom{n}{k}^2 = \binom{2n}{n}$, but there are no simple formulas for $\sum_k \binom{n}{k}^r$ for $r \geq 2$. 
In contrast, Stanley showed that $\sum_k \stern{n}{k}\!^r$ obeys a homogeneous linear recurrence for all $r$.

We repeat Stanley's argument. For any function $f(x,y)$ from $\RR^2 \to \RR$, put 
\[ S_n(f) = \sum_k f\left( \stern{n}{k}, \stern{n}{k+1} \right) . \] 
Define  $(\sigma^{\ast} f)(x,y) = f(x+y,y)$ and $(\tau^{\ast} f)(x,y) = f(x,x+y)$. Then we have
\[\begin{array}{c}
S_n(f) =  \sum_i f\left( \stern{n}{2i}, \stern{n}{2i+1} \right) + \sum_j f\left( \stern{n}{2j+1}, \stern{n}{2j+2} \right) =  \\[1 em]
 \sum_i f\left( \stern{n-1}{i}, \stern{n-1}{i} +  \stern{n-1}{i+1} \right) + \sum_j f\left(\stern{n-1}{j} +  \stern{n-1}{j+1},  \stern{n-1}{j+1} \right)  \\[1 em]
 = S_{n-1}(\sigma^{\ast} f) + S_{n-1}(\tau^{\ast} f).
 \end{array} \]
 
If $f$ is a homogenous polynomial of degree $r$, then $\sigma^{\ast} f$ and $\tau^{\ast} f$ are also such polynomials. 
Thus $f \mapsto \sigma^{\ast} f + \tau^{\ast} f$ is a linear endomorphism of the vector space of degree $r$ polynomials, given by some $(r+1) \times (r+1)$ matrix $\Phi$. And it follows that, for any degree $r$ polynomial $f$, the sequence $S_n(f)$ obeys the same linear recurrence as the sequence of matrix powers $\Phi^n$.

\begin{eg}
We consider cubic polynomials. A basis for the space of homogenous cubics in two variables is $x^3$, $x^2 y$, $x y^2$, $y^3$, and we have
\[ \begin{array}{c@{}c@{}ccc@{}c@{}ccc@{}c@{}c@{}c@{}c@{}c@{}c}
\sigma^{\ast}(x^3) &+& \tau^{\ast}(x^3) &=& (x+y)^3 &+& x^3 &=& 2 x^3 &+& 3 x^2 y &+& 3 x y^2 &+& y^3 \\[1 em]
\sigma^{\ast}(x^2 y) &+& \tau^{\ast}(x^2 y) &=& (x+y)^2 y &+& x^2 (x+y) &=& x^3 &+& 2 x^2 y &+& 2 x y^2 &+& y^3  \\[1 em]
\sigma^{\ast}(x y^2) &+& \tau^{\ast}(x y^2) &=& (x+y) y^2 &+& x (x+y)^2 &=&  x^3 &+& 2x^2 y &+& 2 x y^2 &+& 2 y^3 \\[1 em]
\sigma^{\ast}(y^3) &+& \tau^{\ast}(y^4) &=& y^3 &+& (x+y)^3 &=& x^3 &+& 3 x^2 y &+& 3 x y^2 &+& 3y^3  \\
\end{array} .\]
So $f \mapsto \sigma^{\ast}(f) + \tau^{\ast}(f)$ acts on this four dimensional vector space by the matrix
\[ \Phi= 
\begin{bmatrix}
2 & 1 & 1 & 1 \\
3 & 2 & 2 & 3 \\
3 & 2 & 2 & 3 \\ 
1 & 1 & 1 & 2 \\
\end{bmatrix} .\]
\end{eg}

We say that a sequence $T_n$ obeys a homogenous linear recurrence of length $\ell$ if there are constants $a_1$, $a_2$, \dots, $a_{\ell}$ such that $T_n = \sum_{j=1}^n a_j T_{n-j}$. 
Since the matrix $\Phi$ in the above example is $4 \times 4$, the matrix powers $\Phi^n$ obey a homogenous linear recurrence of length $4$, and thus the $S_n(f)$ do likewise for any homogenous cubic polynomial $f$.
More generally, if $f$ is a homogenous polynomial of degree $r$, this argument shows that the sequence $S_n(f)$ will obey a recurrence of length $r+1$.

This is an elegant argument, but it falls far short of the truth. Computation shows that
\[
S_{\bullet}(x^3) = S_{\bullet}(y^3) =  1,\ 3,\ 21,\ 147,\ \cdots \qquad S_{\bullet}(x^2 y)  = S_{\bullet}(x y^2) = 0,\ 2,\ 14,\ 98,\ \cdots .
\]
We actually have the length $1$ recurrence $S_n(f) = 7 S_{n-1}(f)$ (for $n \geq 2$)!

Why do these recurrences have length $1$ rather than $4$? There are two reasons. The first is that the operator $\Phi$ commutes with the map $f(x,y) \mapsto f(y,x)$, and we have $S_n(f(x,y)) = S_n(f(y,x))$. So $\Phi$ passes to the $2$-dimensional quotient vector space where we identify $f(x,y)$ and $f(y,x)$. Writing $[f]$ for the class of a polynomial in this quotient, we have $\Phi\left( [x^3] \right) = 3 [x^3 ] + 6 [x^2 y]$ and $\Phi\left( [x^2 y ] \right) = 2 [ x^3 ] + 4 [ x^2 y ]$, so $\Phi$ acts by the matrix $\Phi_{\sym} = \begin{sbm} 3&2 \\ 6&4 \end{sbm}$.  Thus $S_n(f)$ obeys a recurrence of length $2$. In general, this argument shows that, for $f$ a homogenous polynomial of degree $r$, the sequence $S_n(f)$ will obey a recurrence of length $\lceil \tfrac{r+1}{2} \rceil$. 

The second reason is that the characteristic polynomial of $\Phi_{\sym}$ is $x^2 - 7x$, so we have $\Phi_{\sym}^n = 7 \Phi_{\sym}^{n-1}$ for $n \geq 2$, and thus $S_n(f) = 7 S_{n-1}(f)$. 
The $0$-eigenvalue makes the recurrence shorter!

Stanley conjectured that the analogous operators $\Phi_{\sym}$ had repeated eigenvalues of $0$ for $r$ odd, and of $\pm 1$ for $r$ even.
The latter phenomenon makes the minimal polynomial of $\Phi_{\sym}$ of lower degree than the characteristic polynomial, and thus also shortens the length of the recurrence satisfied by $S_n(f)$. After introducing some general notation, we state Stanley's conjectures, and prove them.

Fix a nonnegative integer $r$. Let $V$ be the vector space of degree $r$ homogeneous polynomials in $x$ and $y$.
For a  $2 \times 2$ integer matrix $\gamma = \begin{sbm} a&b \\ c&d \end{sbm}$, we define
\[ (\gamma^{\ast} f) (x,y) = f(ax + cy,\ bx+dy) .\]
We put $\sigma = \begin{sbm} 1&0 \\ 1&1 \end{sbm}$ and $\tau =  \begin{sbm} 1&1 \\ 0&1 \end{sbm}$, so this is consistent with our previous notation, and we note that $(\alpha \beta)^{\ast} (f) = \beta^{\ast} \alpha^{\ast} f$. 
We put $\Phi = \sigma^{\ast} + \tau^{\ast}$, a linear map $V \to V$.
 
Let $V_{\sym}$ be the quotient of $V$ by $f(x,y) \equiv f(y,x)$. 
The map $\Phi$ passes to the quotient $V_{\sym}$, and we denote this quotient map by $\Phi_{\sym}$. 
Write $m(M, \lambda)$ for the multiplicity of $\lambda$ as an eigenvalue of the matrix $M$. 

Our main theorem involves many periodic functions of $r$, where $r$ is restricted to either even or odd numbers.
We will write $[a_1, a_2, \ldots, a_p]_{r}$ for the function of such $r$ which is periodic modulo $2p$ with values $a_1$, $a_2$, \dots, $a_p$, $a_1$, $a_2$, \dots, $a_p$, \dots, where we take $a_1$ to occur for $r \equiv 0 \bmod 2p$ if $r$ is restricted to even classes and we take $a_1$ to occur for $r \equiv 1 \bmod 2p$ if $r$ is restricted to even classes.

\begin{MT} 
The operators $\Phi$ and $\Phi_{\sym}$ are diagonalizable over $\RR$.
For $r$ odd, we have
\[ \textstyle{m(\Phi, 0) \geq \tfrac{r}{3} + \left[ -\tfrac{1}{3}, 1, \tfrac{1}{3} \right]_r \qquad m(\Phi_{\sym}, 0) \geq  \tfrac{r}{6} + \left[ -\tfrac{1}{6}, \tfrac{1}{2}, \tfrac{1}{6} \right]_r} .\]
For $r$ even, we have
\[ \begin{array}{l@{\quad}l}
 m(\Phi, 1) \geq \tfrac{r}{6} + \left[-1, \tfrac{2}{3}, \tfrac{1}{3}, 0, -\tfrac{1}{3}, \tfrac{4}{3} \right]_r &  m(\Phi_{\sym}, 1) \geq \tfrac{r}{12} + \left[-1, -\tfrac{1}{6}, -\tfrac{1}{3}, -\tfrac{1}{2}, -\tfrac{2}{3}, \tfrac{1}{6} \right]_r \\[1 em]
 m(\Phi, -1) \geq \tfrac{r}{6} + \left[ 0,-\tfrac{1}{3},\tfrac{4}{3},-1,\tfrac{2}{3}, \tfrac{1}{3} \right]_r & m(\Phi_{\sym}, -1) \geq \tfrac{r}{12} +\left[ 0,-\tfrac{1}{6},\tfrac{2}{3},-\tfrac{1}{2},\tfrac{1}{3},\tfrac{1}{6}\right]_r \\[1 em]
 m(\Phi,1) + m(\Phi,-1) \geq \tfrac{r}{3} + \left[-1, \tfrac{1}{3}, \tfrac{5}{3} \right]_r &   m(\Phi_{\sym},1) + m(\Phi_{\sym},-1) \geq \tfrac{r}{6} + \left[ -1, -\tfrac{1}{3}, \tfrac{1}{3} \right]_r . \\[1 em]
 \end{array} \]
 \end{MT}
  The parts of this theorem related to $\Phi_{\sym}$ were conjectured by Stanley, except that Stanley only conjectured diagonalizability over $\CC$.
  Stanley did not consider $\Phi$, but we think it is natural to consider this operator both because it is easier to compute with and because it would occur naturally if one considered the Stern recurrence with asymmetric initial conditions. 
  Computations show that the lower bounds are equalities for all $0<r\leq 100$.
  
The Main Theorem has the following implications for shortening the recurrence obeyed by $S_n(f)$.
 
 \begin{MC}
For $r$ odd, the sequence $S_n(f)$ obeys a homogenous linear recurrence of length $\leq \tfrac{r}{3} + [\tfrac{2}{3}, 0, \tfrac{1}{3}]_r$ for $n$ large. 
For $r$ even, the sequence $S_n(f)$ obeys a homogenous linear recurrence of length $\leq \tfrac{r}{3} + [4, \tfrac{10}{3}, \tfrac{8}{3}]_r$ for $n$ large. This can be shortened to length  $\leq \tfrac{r}{3} + [2, \tfrac{4}{3}, \tfrac{2}{3}]_r$ if we allow recurrences of the form $S_n(f) = \sum_{j=1}^{\ell} a_j S_{n-j}(f) + b + c (-1)^n$. 
\end{MC}

\section{Diagonalizability of $\Phi$ and $\Phi_{\sym}$}

In this section, we will show that $\Phi$ on $V$ is diagonalizable over $\RR$. Since $\Phi_{\sym}$ is the action of $\Phi$ on the quotient $V_{\sym}$, this also establishes diagonalizability of $\Phi_{\sym}$. Working explicitly in the basis $x^b y^{r-b}$ for $V$, we have
\[ \sigma^{\ast}(x^b y^{r-b}) = \sum_a \binom{b}{a} x^a y^{r-a} \ \mbox{and}\  \tau^{\ast}(x^b y^{r-b}) = \sum_a \binom{r-b}{r-a} x^a y^{r-a} \]
so, in this basis, the entries of the matrix $\Phi$ are given by
\[ \Phi_{ab} = \binom{b}{a} + \binom{r-b}{r-a} . \]
\begin{eg}
Returning to the example of cubic polynomials, we have
\[ \Phi= 
\begin{bmatrix}
2 & 1 & 1 & 1 \\
3 & 2 & 2 & 3 \\
3 & 2 & 2 & 3 \\ 
1 & 1 & 1 & 2 \\
\end{bmatrix} =
\begin{bmatrix}
1 &  &  &  \\
3 & 1 &  &  \\
3 & 2 & 1 &  \\ 
1 & 1 & 1 & 1 \\
\end{bmatrix} + 
\begin{bmatrix}
1 & 1 & 1 & 1 \\
 & 1 & 2 & 3 \\
 & & 1 & 3 \\
 & & & 1 \\ \end{bmatrix} . \] 
\end{eg}

Conjugate $\Phi$ by the diagonal matrix $D$ with entries $D_{kk} = \sqrt{k! (r-k)!}$.
For $a \neq b$, we have
\[ \begin{array}{rcl}
(D \Phi D^{-1})_{ab} &=& \sqrt{\tfrac{a! (r-a)!}{b! (r-b)!}} \cdot \begin{cases} \binom{b}{a} & \mbox{for}\ a<b \\ \binom{r-b}{r-a} &\mbox{for}\  a>b \\ \end{cases} \\[0.4 cm]
&=& \tfrac{1}{|a-b|!} \sqrt{\tfrac{\max(a,b)! \max(r-a, r-b)!}{\min(a,b)! \min(r-a,r-b)!}} . \\
\end{array} \]
The last expression is symmetric in $a$ and $b$, so $D\Phi D^{-1}$ is a symmetric matrix and thus diagonalizable over $\RR$. We conclude that $\Phi$ is likewise diagonalizable over $\RR$.

\section{The case of $r$ odd}

Throughout this section, we take $r$ to be odd.
Set
\[ \rho = \sigma \tau^{-1} = \begin{bmatrix} 1 & -1 \\ 1&0 \end{bmatrix} . \]
Let $\omega$ be a primitive $6$-th root of unity. Then $\rho$ acts on $\CC^2$ with eigenvalues $\omega$ and $\omega^{-1}$.
Thus $\rho^{\ast}$ acts on $V \otimes \CC$ with eigenvalues $\omega^a (\omega^{-1})^{r-a}$ for $0 \leq a \leq r$ for $0 \leq a \leq r$; the corresponding eigenvectors are $(x+\omega y)^a (\omega x + y)^{r-a}$.
In particular, the $(-1)$-eigenspace of $\rho^{\ast}$ on $V$ corresponds to those values of $a$ for which $a \equiv r-a+3 \bmod 6$ and has dimension equal to the number of such $a$, which is $\tfrac{r}{3} + [ -\tfrac{1}{3}, 1, \tfrac{1}{3}]_r$.

Let $W$ be this $(-1)$-eigenspace of $\rho^{\ast}$.
For $g \in W$, we have $\rho^{\ast} g = - g$, so $\sigma^{\ast} g = - \tau^{\ast} g$. Then $\Phi(g) = \sigma^{\ast} g+\tau^{\ast} g = 0$.
So $\Phi$ acts by $0$ on $W$, and thus $m(\Phi,0)  \geq \tfrac{r}{3} + [ -\tfrac{1}{3}, 1, \tfrac{1}{3}]_r$

We now consider the corresponding computation for $\Phi_{\sym}$. 
Let $W_{\sym}$ be the $(-1)$-eigenspace of $\Phi_{\sym}$ on $V_{\sym}$. 
In $\Phi_{\sym}$, the eigenvectors $(x+\omega y)^a (\omega x + y)^{r-a}$ and $(x+\omega y)^{r-a} (\omega x + y)^{a}$ are identified, and no other relations are imposed. So the dimension of $W_{\sym}$ is the number of $a$ for which $a \equiv r-a+3 \bmod 6$, modulo the relation that $a$ and $r-a$ are equivalent; the number of such $a$ is  $\tfrac{r}{6} + \left[ -\tfrac{1}{6}, \tfrac{1}{2}, \tfrac{1}{6} \right]_r$.
We have proved the multiplicity bound for $m(\Phi_{\sym},0)$.

\section{The case of $r$ even}

Throughout this section, we take $r$ to be even.
We continue to use the notation $\rho = \sigma \tau^{-1} = \begin{sbm} 1 & -1 \\ 1&0 \end{sbm}$ and also introduce 
\[ \iota =  \sigma  \tau^{-1}  \sigma = \begin{bmatrix} 0&-1 \\ 1&0 \end{bmatrix} .\]
The eigenvalues of $\rho$ on $\CC^2$ are the primitive $6$-th roots of unity and the eigenvalues of $\iota$ are $\pm i$.
We deduce that, on $V = \mathrm{Sym}^r \RR^2$, we have $(\rho^{\ast})^3 = (\iota^{\ast})^2 = \mathrm{Id}$. 

Let $X$ be the subspace of $V$ where $(\rho^{\ast})^2 + (\rho^{\ast}) + 1 =0$, and let $Y^+$ and $Y^-$ be the subspaces of $V$ where $\iota^{\ast}$ has eigenvalues $1$ and $-1$ respectively. We have
\[ \dim X =  \tfrac{2r}{3} + [0, \tfrac{2}{3}, \tfrac{4}{3}]_r  \qquad \dim Y^+ = \tfrac{r}{2} + [1,0]_r \qquad \dim Y^- = \tfrac{r}{2} +[0,1]_r . \]
We deduce that

\[ \begin{array}{lclcl}
\dim (X \cap Y^+) &\geq&\tfrac{2r}{3} + [0, \tfrac{2}{3}, \tfrac{4}{3}]_r + \tfrac{r}{2} + [1,0]_r - (r+1) &=& \tfrac{r}{6} + [0, -\tfrac{1}{3}, \tfrac{4}{3}, -1, \tfrac{2}{3}, \tfrac{1}{3}]_r
 \\[1 em]
\dim (X \cap Y^-) &\geq&\tfrac{2r}{3} + [0, \tfrac{2}{3}, \tfrac{4}{3}]_r + \tfrac{r}{2} + [0,1]_r - (r+1) &=& \tfrac{r}{6} + [-1, \tfrac{2}{3}, \tfrac{1}{3}, 0, -\tfrac{1}{3}, \tfrac{4}{3}]_r.
\end{array}  \]
Adding these formulas,
\[ \textstyle{\dim (X \cap Y^+) + \dim (X \cap Y^-)  \geq  \tfrac{r}{3} + [ -1, \tfrac{1}{3}, \tfrac{5}{3} ]_r } . \]
 
For $f$ in $X$, we have
\[  \Phi(f) = \sigma^{\ast}(f) + \tau^{\ast}(f) = \tau^{\ast} \left( \rho^{\ast} + 1 \right) (f) =  - \tau^{\ast} (\rho^{\ast})^2 (f) = - \sigma^{\ast} (\tau^{-1})^{\ast} \sigma^{\ast}(f) = - \iota^{\ast} (f). \]
So, if $f \in X \cap Y^{\pm}$, we have $\Phi(f) = \mp f$. 
This proves the claimed bounds for $m(\Phi, \pm 1)$.

We now repeat the analysis for $\Phi_{\sym}$. 
We define $X_{\sym}$ and $Y^{\pm}_{\sym}$ analogously. 
A basis for $X \otimes \CC$ is those polynomials $(x+ \omega y)^a (\omega x + y)^{r-a}$ for which $2a \not \equiv r \bmod 3$. In $X_{\sym}$ the polynomials $(x+ \omega y)^a (\omega x + y)^{r-a}$ and $(x+ \omega y)^{r-a} (\omega x + y)^{r}$ are identified, which cuts the dimension exactly in half. (Note that $(x+ \omega y)^{r/2} (\omega x + y)^{r/2}$ has eigenvalue $1$, so it doesn't contribute to $X_{\sym}$.)
So 
\[ \dim X_{\sym} =  \tfrac{r}{3} + [0, \tfrac{1}{3}, \tfrac{2}{3}]_r . \]
Similar but simpler analysis shows that 
\[ \dim Y^+_{\sym} = \tfrac{r}{4}+[1,\tfrac{1}{2}]_r  \qquad \dim Y^-_{\sym} = \tfrac{r}{4}+[0,\tfrac{1}{2}]_r . \]

We have
\[ \begin{array}{lclcl}
\dim (X_{\sym} \cap Y_{\sym}^+) &\geq&\tfrac{r}{3} + [0, \tfrac{1}{3}, \tfrac{2}{3}]_r+ \tfrac{r}{4}+[1,\tfrac{1}{2}]_r -(\tfrac{r}{2}+1)&=&  \tfrac{r}{12}+ [0, -\tfrac{1}{6}, \tfrac{2}{3}, -\tfrac{1}{2}, \tfrac{1}{3}, \tfrac{1}{6}]_r \\[1 em]
\dim (X_{\sym} \cap Y_{\sym}^-) &\geq&\tfrac{r}{3} + [0, \tfrac{1}{3}, \tfrac{2}{3}]_r+ \tfrac{r}{4}+[0,\tfrac{1}{2}]_r -(\tfrac{r}{2}+1)&=&  \tfrac{r}{12}+ [-1, -\tfrac{1}{6}, -\tfrac{1}{3}, -\tfrac{1}{2}, -\tfrac{2}{3}, \tfrac{1}{6}]_r \\[1 em]
\end{array}  \]
Adding these,
\[ m(\Phi_{\sym},1) + m(\Phi_{\sym},-1) \geq \tfrac{r}{6} + \left[ -1, -\tfrac{1}{3}, \tfrac{1}{3} \right]_r  . \]
We have now deduced all claims about $m(\Phi_{\sym}, \pm 1)$.

\section{Consequences for the length of recurrences}

First, we consider the case of $r$ odd.
The matrix $\Phi_{\sym}$ is square of size $\tfrac{r+1}{2}$; let the characteristic polynomial of $\Phi_{\sym}$ be $x^{(r+1)/2} - \sum_{j=1}^{(r+1)/2} a_j x^{(r+1)/2-j}$.
So we have
\[ \Phi_{\sym}^n = \sum_{j=1}^{(r+1)/2} a_j \Phi_{\sym}^{n-j}  \]
But the last $m(\Phi_{\sym},0)$ terms of the sum are $0$, so the sum actually only runs up to $j = \tfrac{r+1}{2} - m(\Phi_{\sym},0) \geq \tfrac{r}{3} + [ \tfrac{2}{3}, 0, \tfrac{1}{3}]_r$.
For any degree $r$ polynomial $f$, the sequence $S_n(f)$ is a linear function of $\Phi_{\sym}^n$, so $S_n(f)$ obeys a homogenous linear recurrence of the same length.

We now consider the case of $r$ odd. Write the characteristic polynomial of $\Phi_{\sym}$ in the form $g(x) (x-1)^{m(\Phi_{\sym},1)} (x+1)^{m(\Phi_{\sym},1)}$. 
Then $\Phi_{\sym}$ obeys the polynomial $g(x) (x-1) (x+1)$ with degree $\tfrac{r}{2} + 1 - \left( m(\Phi_{\sym},1) + m(\Phi_{\sym},-1) -2 \right) =  \tfrac{r}{3} + [4, \tfrac{10}{3}, \tfrac{8}{3}]_r$. Hence $\Phi_{\sym}^n$ obeys a homogenous linear recurrence of this length, and $S_n(f)$ does as well.

Moreover, change bases to put $\Phi_{\sym}$ into block form $\begin{sbm} M & 0 & 0 \\ 0 & \mathrm{Id} & 0 \\ 0 & 0 & -\mathrm{Id} \\ \end{sbm}$, so $g(M) = 0$ and let $g(x) = \sum g_j x^{\deg(g) - j}$. 
Then $\sum g_j \Phi_{\sym}^{n-j}$ will have block form $\begin{sbm} 0 & 0 & 0 \\ 0 & p \mathrm{Id} & 0 \\ 0 & 0 & q (-1)^n \\ \end{sbm}$ and hence $\sum g_j S_{n-j}(f)$ will be of the form $b + c(-1)^n$ for some constants $b$ and $c$. 
So $S_n(f)$ obeys  a recurrence of the form $S_n(f) = \sum_{j=1}^{\ell} a_j S_{n-j}(f) + b + c (-1)^n$ for $\ell = \deg(g)$.

\section{Acknowledgments} 
The author was partially supported by NSF grant DMS-1600223. 
The author learned of this conjecture from Richard Stanley's talk at ``Combinatorics and beyond", a conference in honor of Sergey Fomin, and thanks Professor Stanley for bringing this delightful construction to his attention.

\end{document}